\documentclass[12pt,letterpaper]{article}
\usepackage[latin1]{inputenc}
\usepackage{amsmath}
\usepackage{amsfonts}
\usepackage{amssymb}
\usepackage{amsthm}

\usepackage{ragged2e}

\title{The Nuisance Principle in Infinite Settings }
 \author{Sean C. Ebels-Duggan}
\date{Submitted 1 July 2015}
\begin{document}

\maketitle

\paragraph{Disclaimer}  This is the pre-peer-reviewed version of the following article:
\begin{quote}
Sean C. Ebels-Duggan.  The Nuisance Principle in Infinite Settings, \emph{Thought: A Journal of Philosophy}~4(4), December~2015, pp.~263--268.
\end{quote} 
which has been published in final form at 
\begin{quote}
http://onlinelibrary.wiley.com/doi/10.1002/tht3.186/abstract
\end{quote}

\paragraph{Note}  The final version has additional remarks and corollaries; namely these: First corollary: SOL + Nuisance Principle (NP) prove that there is no pairing injection of the universe. Second corollary: SOL+ HP + NP prove that the universe is uncountable. Remark: Even if NP does not imply the finitude of the universe, it is still deductively non-conservative.)
%\RaggedRight

\section*{Submitted Article}

\begin{abstract}
    Neo-Fregeans have been troubled by the Nuisance Principle (NP), an abstraction principle that is consistent but not jointly (second-order) satisfiable with the favored abstraction principle HP. We show that logically this situation persists if one looks at joint (second-order) consistency rather than satisfiability: under a modest assumption about infinite concepts, NP is also inconsistent with HP. 
\end{abstract}

The so-called Nuisance Principle (NP) is the paradigm example of an abstraction principle that is individually satisfiable in second-order logic (with full comprehension), but is not jointly satisfiable with the neo-logicist's favored abstraction principle, HP. This is thought to cause trouble for neo-logicists. Some abstraction principles, when added to second-order logic, allow the recovery of certain mathematical content. For example, HP allows one to recover second-order Peano Arithmetic. If abstraction principles have epistemic status near enough to logic, then so do the recovered mathematics. The principle NP is troublesome because initially it seems to have epistemic status like HP, but it is hard to see how near-logical principles could be so incompatible. 

But are NP and HP jointly \emph{consistent}?  The further question arises because satisfiability (having a \emph{standard} model) and consistency (not proving a contradiction) %, and by completeness, having \emph{some} model) 
are not the same in second-order logic.  The question was partially answered %by Walsh 
in \cite[21--22]{WalshED2015};  the present note moves us further, but not fully, towards a complete answer.

The principle HP, attributed loosely to Hume by Frege \cite[$\S$ 63]{Frege1980}, states that the Number of $F$s (denoted $\#F$) is identical to the Number of $G$s ($\#G $) just if there is a bijection from the $F$s to the $G$s---that is, a function associating all of the objects falling under $F$ with all of the objects falling under $G$, such that no two objects falling under $F$ are associated with the same object falling under $G$.  In second-order logic the existence of such a bijection can be represented, and is demonstrably an equivalence relation.  Thus HP can be represented by:
\[ (\forall F)(\forall G)(\#F = \#G \leftrightarrow F \approx G)\]
where `$\approx$' is shorthand for the second-order formula asserting the existence of a bijection.  

The Nuisance Principle is a simplification due to Crispin Wright \cite{Wright1997aa} of a principle introduced by George Boolos \cite{Boolos1990}.  One can express in second-order logic the following equivalence relation:  %$N(X,Y)$ holds if and only if 
\begin{quote} $N(F,G)$ \emph{iff} there are finitely many objects falling under $F$ but not $G$, and finitely many falling under $G$ but not $F$
\end{quote}
The Nuisance principle is then the claim that the \emph{Nuisance} of $F$ (denoted $\ddag F$) is identical to the Nuisance of $G$ ($\ddag G$) if and only if $N$ holds of $F$ and $G$.  Using our abbreviation $N(F, G)$ we can represent this in second-order logic by
\[ (\forall F)(\forall G)(\ddag F = \ddag G \leftrightarrow N(F,G) )\]
Notice that NP and HP are both abstraction principles in virtue of having the same form:  equality between objects on the left, an equivalence relation between concepts on the right.

That NP and HP are jointly unsatisfiable can be seen by deploying features of cardinal numbers in set theory to show that the former is satisfiable only if there are finitely many objects.\footnote{For such a proof that NP is unsatisfiable, see \cite{Antonelli2010aa}.}  Since HP proves there are infinitely many objects, the two are not jointly satisfiable.  But one cannot adapt this proof to a deductive setting. The complicating factor is that outside of standard models, being ``infinite'' can mean many things.  Typically, concepts are ``infinite'' if they are \emph{Dedekind infinite}:  there is a function from \emph{all} the objects falling under the concept to \emph{not all} of the objects falling under the concept, such that no two objects are sent to the same object by that function.  (That is to say, there is an injection from the concept to a proper subconcept of itself.)  In standard models of second-order logic, Dedekind infinite concepts behave like infinite sets behave in set theory.  But this isn't guaranteed in non-standard models (and this is why in this note we use ``concepts'' rather than ``sets'' to indicate the semantic correlate to second-order variables).

In this note we show that NP is inconsistent with the Dedekind infinity of the universe in the presence of a natural and relatively modest \emph{strengthening} of the assertion that the universe is Dedekind infinite.  Such a strengthening is a conditional describing the behavior of Dedekind infinite concepts.  %whose antecedent is the second-order formula expressing that the universe is Dedekind infinite.  

This is a significant improvement over what was shown in \cite{WalshED2015}.  The proof in that paper used two versions of the Axiom of Choice:  a global well-ordering {\tt GC} to get Dedekind infinite concepts to behave like infinite sets, and a uniform means of selecting representatives for each equivalence class, {\tt AC}.   So what was shown in that paper is that, 
%What %Walsh showed 
%was shown in \cite{WalshED2015} is that, 
if one's second-order logic includes these versions of the Axiom of Choice, then NP is not consistent with the universe being Dedekind infinite.  Thus HP and NP are jointly inconsistent, as in the proof that they are unsatisfiable.  

Our improvement is that we can obtain this result by appeal to an ostensibly weaker principle.  
The principle in question is the following strengthening of infinity:
\begin{description}
\item[(Pairing)]  If the universe is Dedekind infinite, then there is a binary function $f$ defined on all pairs of objects such that for any $z$ and any $x,y,x',$ and $y'$, if $z=f(x,y)$ and $z=f(x', y')$, then $x=x'$ and $y= y'$. %\footnote{We help ourselves to functions in second-order logic, but note that the function $f$ described above can be replaced by a relation $R$ that is the graph of the described $f$.
%The versions of the Axiom of Choice used in \cite{WalshED2015}were a global well-ordering of the universe (there identified as {\tt GC}), and a means of selecting a representative from each equivalence class (identified as {\tt AC}). It is also worth noting with \cite{WalshED2015} that Pairing is a consequence of {\tt GC}, and indeed that in {\tt ZF}, a version of Pairing implies the (set theoretic) Axiom of Choice (see \cite[Theorem 11.7]{Jech1973aa}).  But equivalence in {\tt ZF} is not the same as equivalence in second-order logic;  for this reason these principles must be treated separately.}
%		\item $R$ is the graph of a function $f$ on pairs:  For any objects $x$ and $y$ there is a unique  $z$ such that $f(x,y)=z$ 
%$f$ is injective:  
%		\item $f$ is surjective:  For any $z$, there are $x$ and $y$ such that $f(x,y)=z$.
 %\end{enumerate}
\end{description}
%The principle Pairing is an axiom of infinity:  fix any object $a$, then the one-place function $f(a,y)$ witnesses that the universe is Dedekind infinite.  %In what follows, for objects $x$ and $y$ we will denote that unique $z$ such that $R(x,y,z)$ by ``$\langle x, y\rangle$''.
In other words, if the universe is Dedekind infinite, then there is an injection from pairs of objects into the universe.\footnote{It is worth reiterating the remark of \cite{WalshED2015} that Pairing is a consequence of {\tt GC}.  It is also worth the separate remark that in {\tt ZF} set theory, a version of Pairing implies the (set theoretic) Axiom of Choice (see \cite[Theorem 11.7]{Jech1973aa}). Because equivalence in {\tt ZF} is not the same as equivalence in second-order logic, we here treat these these principles as distinct.}  In effect, this strengthening says that universe-sized concepts can be broken-up into universe-many disjoint subconcepts, each of universe-size.  

We now sketch a deductive argument showing that NP and Pairing are inconsistent with the assertion that the universe is Dedekind infinite.  For if the universe is Dedekind infinite and Pairing holds,  %then for each %object $a$ we can associate a concept $S[a]$ as follows:
%\[ z \; \textrm{falls under}\; S[a] \leftrightarrow \; \textrm{there is a $y$ such that}\; f(a,y)=z \]
%In other words, $S[x]$ is the image of the universe under the one-place function $f(a,y)$.
%
%Now using $S$ 
we can associate a Dedekind infinite concept with each concept, \emph{whether the latter is finite or infinite}.  For given a concept $F$, let $U[F]$ be defined by 
\begin{multline*} z \; \textrm{falls under}\; U[F] \leftrightarrow \\
\textrm{there is an $x$ falling under $F$, and a $y$ such that}\; f(x,y)=z \end{multline*}
In other words, $U[F]$ is the image of $F$ when projected (on the right) with the universal set $V$ (on the left):  $U[F] = f(F,V)$.

Now we show that if concepts $F$ and $G$ are extensionally distinct (if something falls under one that doesn't fall under the other), then the equivalence relation $N$, described above, does not hold between $U[F]$ and $U[G]$.  For if $a$ falls under $F$ but not $G$, then by fixing $a$ we obtain a one-to-one map $f(a, y)$ from the universal set into $U[F]-U[G]$, the part of $U[F]$ that does not overlap $U[G]$.  Since the universe is Dedekind infinite, by Pairing, so is $U[F]-U[G]$.  An identical argument can be made for any element falling under $G$ but not $F$.  Thus $N$ does not obtain between $U[F]$ and $U[G]$.  Clearly, of course, if $F$ and $G$ are not extensionally distinct, then $N(U[F], U[G])$, since $U[A]$ and $U[B]$ will not be extensionally distinct either.  In other words, 
\[ (\forall F)(\forall G)(N(U[F], U[G]) \leftrightarrow (\forall x)(Fx \leftrightarrow Gx)) \]

Suppose now that NP obtains;  we then have %using the the operator $\ddag$ to obtain that for any concepts $F$ and $G$ (finite or not):
%\[ \ddag U[X] = \ddag U[Y] \leftrightarrow N(U[A], U[B]) \leftrightarrow (\forall x)(Xx \leftrightarrow Yx) \]
%That is to say:
\[ (\forall F)(\forall G)(\ddag U[F] = \ddag U[G] \leftrightarrow (\forall x)(Fx \leftrightarrow Gx)) \]
One can then rehearse the argument of Russell's paradox:  where $R$ is the concept defined by 
\begin{multline*} y \; \textrm{falls under}\; R \leftrightarrow 
\textrm{there is a concept $Y$ such that $y = \ddag U[Y]$} \\
\textrm{and $y$ does not fall under $Y$},\end{multline*}
the usual argument shows that $\ddag U[R]$ falls under $R$ if and only if it doesn't.  So Pairing and NP imply that the universe is not Dedekind infinite.  Thus in the presence of Pairing, HP and NP are not jointly consistent.\footnote{Acknowledgments removed for blind review.
}

\bibliography{Nuis-bib}
\bibliographystyle{plain}

\end{document}